\documentclass[12pt]{article}

\usepackage{latexsym,amssymb,amsmath,theorem,ifthen,epsfig,diagrams,eepic}
\usepackage{makeidx,manfnt,epsfig}
\usepackage{ulem,stmaryrd}

\newboolean{LONG}\newboolean{EXCLUDE}
\setboolean{EXCLUDE}{true}

\def\dotuline{\bgroup 
  \ifdim\ULdepth=\maxdimen  
   \settodepth\ULdepth{(j}\advance\ULdepth.4pt\fi
  \markoverwith{\begingroup
  \advance\ULdepth0.08ex
  \lower\ULdepth\hbox{\kern.15em .\kern.1em}%
  \endgroup}\ULon}

\def\dashuline{\bgroup 
  \ifdim\ULdepth=\maxdimen  
   \settodepth\ULdepth{(j}\advance\ULdepth.4pt\fi
  \markoverwith{\kern.15em
  \vtop{\kern\ULdepth \hrule width .3em}%
  \kern.15em}\ULon}

\title{A remark on the invariance of dimension}
\author{Michael M\"uger \\ Radboud University, Nijmegen, The Netherlands}

\newlength{\dinwidth}
\newlength{\dinmargin}
\def\endexem{\hfill{$\Box$}\medskip}

\def\1#1{{\bf #1}}
\def\2#1{{\cal #1}}
\def\3#1{{\sl #1}}
\def\4#1{{\tt #1}}
\def\5#1{{\sf #1}}
\def\6#1{{\mathfrak #1}}
\def\7#1{{\mathbb #1}}

\def\qed{\hfill{$\blacksquare$}\medskip}

\reversemarginpar

\newcommand{\dist}{\mathrm{dist}}
\newcommand{\diam}{\mathrm{diam}}

\def\prf{\noindent {\it Proof.\ }}
\def\id{\mathrm{id}}

\newcommand{\restr}{\!\upharpoonright\!}
\newcommand{\rarr}{\rightarrow}

\newcommand{\impl}{\Rightarrow}
\newcommand{\ol}{\overline}

\newcommand{\ve}{\varepsilon}

\newtheorem{thm}{Theorem}[section]
\newtheorem{prop}[thm]{Proposition}
\newtheorem{lem}[thm]{Lemma}
\newtheorem{cor}[thm]{Corollary}
\newtheorem{defin}[thm]{Definition}
\newtheorem{defprop}[thm]{Definition/Proposition}
\theorembodyfont{\rmfamily}
\newtheorem{example}[thm]{Example}
\newtheorem{exercise}[thm]{Exercise}
\newtheorem{rema}[thm]{Remark}

\newcommand{\bdefin}{\begin{defin}}
\newcommand{\blemma}{\begin{lem}}
\newcommand{\bprop}{\begin{prop}}
\newcommand{\btheor}{\begin{thm}}
\newcommand{\bcoro}{\begin{cor}}
\newcommand{\bconj}{\begin{conj}}
\newcommand{\bdefprop}{\begin{defprop}}
\newcommand{\bexam}{\begin{example}}
\newcommand{\bexer}{\begin{exercise}}
\newcommand{\edefin}{\end{defin}}
\newcommand{\elemma}{\end{lem}}
\newcommand{\eprop}{\end{prop}}
\newcommand{\etheor}{\end{thm}}
\newcommand{\ecoro}{\end{cor}}
\newcommand{\econj}{\end{conj}}
\newcommand{\brem}{\begin{rema}}
\newcommand{\erem}{\endexem\end{rema}}
\newcommand{\edefprop}{\end{defprop}}
\newcommand{\eexam}{\endexem\end{example}}
\newcommand{\eexer}{\end{exercise}}

\newcommand{\be}{\begin{equation}}
\newcommand{\ee}{\end{equation}}
\newcommand{\bean}{\begin{eqnarray*}}
\newcommand{\eean}{\end{eqnarray*}}

\theorembodyfont{\rmfamily}

\newcommand{\del}{\partial}

\numberwithin{equation}{section}
\numberwithin{thm}{section}

\makeindex

\begin{document}
\maketitle

\abstract{Combining Kulpa's proof of the cubical Sperner lemma and a dimension theoretic idea of van Mill we
  give a very short proof of the invariance of dimension, 
  i.e.\ the statement that cubes $[0,1]^n, [0,1]^m$ are homeomorphic if and only if $n=m$.

This note is adapted from lecture notes for a course on general topology.}

\section{Introduction}
In introductory courses to (general) topology, as also this author has given them many times, it is standard
to show that the connected subsets of $\7R$ are precisely the convex ones (`intervals') and to exhibit the
intermediate value theorem from calculus as an easy consequence. Brouwer's fixed-point theorem for $I=[0,1]$
follows trivially. Connectedness arguments also readily prove that $\7R\not\cong\7R^n, I\not\cong I^n$ if
$n\ge 2$. If the course covers the fundamental group, the non-vanishing of $\pi_1(S^1)$ is used to prove the
fixed-point theorem for the 2-disk and to distinguish $\7R^2$ from $\7R^n$ for $n\ge 3$. But the generalization
of these results to higher dimensions is usually omitted, referring to courses in algebraic topology. Typical
representatives of this approach are \cite{munkres,runde}.

There are, to be sure, proofs of the fixed-point theorem that do not (explicitly) invoke algebraic topology. 
On the one hand, there is a long tradition of proofs that use a combination of calculus and linear algebra,
the nicest perhaps being the one of Lax \cite{lax}. However, these proofs seem vaguely inappropriate
considering that basic topology should be a more elementary subject than calculus. (At this stage, the student
would not be helped by the information that behind such proofs there is de Rham cohomology.)
On the other hand, there is the combinatorial approach via Sperner's lemma
\cite{sperner,KKM}. (Cf.\ \cite[pp.\ 411-417]{engelking} for a nice presentation.) 
But also the combinatorial approach is not entirely satisfactory since it requires introducing simplicial
language and is complicated due its use of barycentric subdivision. Neither of these complaints applies to the
beautiful and elementary proof of the fixed-point theorem published by Kulpa in 1997, cf.\ \cite{kulpa}. But one
still would like to have an accessible proof of the invariance of dimension. 

It is well known that the invariance of dimension can be deduced from Brouwer's fixed-point theorem in at least
two different ways. On the one hand, one can use Borsuk's theory of maps into spheres (cohomotopy) to prove
the `invariance of domain', another result of Brouwer, from which invariance of dimension readily follows. 
This approach is followed, e.g., in Eilenberg's and Steenrod's {\it Foundations of algebraic topology}, in the 
topology books of Kuratowski and of Dugundji, and in Engelking's and Siekluchi's  {\it Topology.\ A geometric
  approach}.  A drawback of this approach is that it invariably requires some use of simplicial techniques
beyond those typically employed in the proof of the fixed-point theorem. In \cite{kulpa2}, Kulpa simplified these
methods somewhat and obtained a proof of the invariance of dimension in about seven pages. However, this is
still more involved than one might wish. 

The second way to deduce the invariance of dimension from the fixed-point theorem relies on dimension theory.
Dimension theory (cf.\ \cite{engelking3} for the most up-to-date account) associates to a space $X$ an element
$\dim(X)\in\{-1,0,1,\ldots,\infty\}$ in a homeomorphism-invariant manner. (There are actually several
competing definitions of dimension.)  
Invariance of the dimension results as soon as one proves $\dim(I^n)=n$. The proof of $\dim(I^n)\le n$ is
easy, but the opposite inequality requires a certain non-separation result, cf.\ Corollary \ref{coro-kulpa}
below. The latter result can be deduced from Brouwer's fixed-point theorem, and this is done in virtually all
expositions of dimension theory. However, setting up enough of dimension theory (for any one of the existing 
definitions) to prove $\dim(I^n)=n$ again requires several pages. In \cite{mill}, van Mill exhibits a very
elegant and efficient alternative approach. He uses yet another notion of dimension, for which proving the
upper bound is very easy, and which has the virtue that Corollary \ref{coro-kulpa} immediately gives the lower
bound. We call this notion the `separation-dimension', cf.\ Definition \ref{sdim}. (For a brief history of the
latter, cf.\ Remark \ref{rem-hist}.)

The purpose and only original contribution of this otherwise expository note is the observation that what
Kulpa `really' proves in \cite{kulpa} is Theorem \ref{theor-kulpa} below, from which Corollary
\ref{coro-kulpa} follows as immediately as does the Poincar\'e-Miranda theorem. This allows to cut out any
reference to the latter (and to Brouwer's theorem) and to give a proof of dimension invariance that easily
fits into four pages and can be explained in one lecture of 90 minutes. This makes it even shorter than
Brouwer's first proof \cite{brouwer-dim}, which was neither self-contained nor easy to read. 
It also bears emphasizing that deducing both the fixed-point theorem and the invariance of dimension from the
higher-dimensional connectedness of the cube asserted by Theorem \ref{theor-kulpa} makes these deductions
entirely analogous to the trivial ones in dimension one.
The price to pay for this efficient approach is that it does not provide a proof of the invariance of domain.

The heart of this note is Section \ref{sec-dim}, but for the reader's convenience, we include an appendix with 
Kulpa's deduction of the theorems of Poincar\'e-Miranda and Brouwer and some corollaries.

\vspace{.5cm}
\noindent{\it Acknowledgments.} The author would like to thank Arnoud van Rooij for suggesting a
simplification of the proof of Proposition \ref{prop-mill}. He also thanks two referees for constructive
comments that led to improvements of the presentation.

\section{The cubical Sperner lemma}
In the entire section, $n\in\7N$ is fixed and $I=[0,1]$. The \underline{faces} of the $n$-cube $I^n$ are given by
\[ I_i^-=\{ x\in I^n\ | \ x_i=0\},\quad\quad I_i^+=\{ x\in I^n\ | \ x_i=1\}. \]
We need some more notations: 

\begin{itemize}
\item For $k\in\7N$, we put $\7Z_k=k^{-1}\7Z=\{ n/k\ | \ n\in\7Z\}$. Clearly $\7Z_k^n\subset\7R^n$.
\item $e_i\in\7Z_k^n$ is the vector whose coordinates are all zero, except the $i$-th, which is $1/k$.
\item $C(k)=I^n\cap \7Z_k^n=\left\{0,\frac{1}{k},\ldots,\frac{k-1}{k},1\right\}^n$. (The
  \underline{combinatorial $n$-cube}.) 
\item $C_i^\pm(k)=I_i^\pm\cap \7Z_k^n$. (The \underline{faces} of the combinatorial $n$-cube $C(k)$.) 
\item $\del C(k)=\cup_i (C_i^+(k)\cup C_i^-(k))$. (The \underline{boundary} of the combinatorial $n$-cube $C(k)$.)
\item A \underline{subcube} of $C(k)$ is a set $C=\{ z_0+\sum_{i=1}^n a_i e_i\ | \ a\in\{0,1\}^n\}\subset C(k)$,
 where $z_0\in C(k)$.
\end{itemize}

Sperner's lemma in its original form \cite{sperner,KKM} concerns simplices. The cubical version stated below
was first proven in \cite{kuhn} and differently in \cite[Lemma 1]{wolsey}. The following proof is the one
given in \cite{kulpa}.  

\bprop \label{prop-sperner}
Let $\varphi: C(k)\rarr\{0,\ldots,n\}$ be a map such that (i) $x\in C_i^-(k)\,\impl\,\varphi(x)<i$ and (ii) 
$x\in C_i^+(k)\,\impl\,\varphi(x)\ne i-1$. Then there is a subcube $C\subset C(k)$ such that
$\varphi(C)=\{0,\ldots,n\}$. 
\eprop

\bdefin An \underline{$n$-simplex} in $\7Z_k^n$ is an ordered set $S=[z_0,\ldots,z_n]\subset\7Z_k^n$ such that
\[ z_1=z_0+e_{\alpha(1)}, \ \ z_2=z_1+e_{\alpha(2)}, \ \ \ldots,\ \ z_n=z_{n-1}+e_{\alpha(n)}, \]
where $\alpha$ is a permutation of $\{1,\ldots,n\}$. 
The subset $F_i(S)=[z_0,\ldots,z_{i-1},z_{i+1},\ldots,z_n]\subset S$, where $i\in\{0,\ldots,n\}$, is called
the $i$-th \underline{face of the $n$-simplex $S$}. 

A finite ordered set $F\subset\7Z_k^n$ is called a \underline{face} if it is a face of some simplex. 
\edefin

Note that the faces $F_i(S)$ are $(n-1)$-simplices in the above sense only if $i=0$ or $i=n$.

\begin{figure}[h]
         \centerline{
           \scalebox{1}{
             \epsfig{file=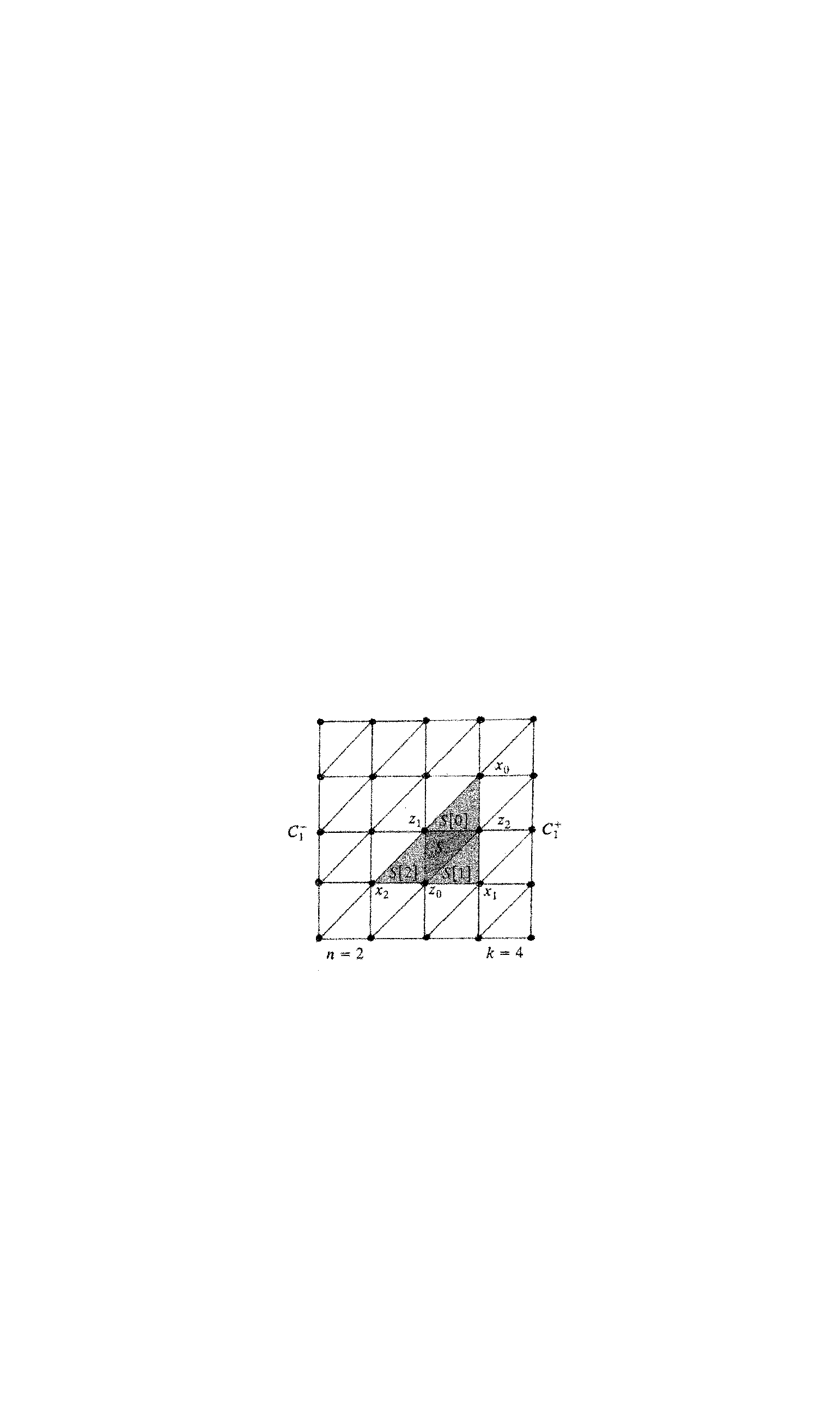, scale=0.9}
             }
           }
         \caption{The neighbors of a simplex in $\7Z^2_4$. From \cite{kulpa} in Amer.\ Math.\ Monthly.}
         \label{fig-kulpa}
 \end{figure}

\blemma \begin{itemize}
\item[(i)] Let $S=[z_0,\ldots,z_n]\subset\7Z_k^n$ be an $n$-simplex. Then for every $i\in\{0,\ldots,n\}$ there
is a unique $n$-simplex $S[i]$, the \underline{$i$-th neighbor} of $S$, such that 
$S\cap S[i]=F_i(S)$.
\item[(ii)] If $S\subset C(k)$ and $i\in\{0,\ldots,n\}$ then $S[i]\subset C(k)$ holds if and only if
  $F_i(S)\not\subset\del C(k)$.  
\end{itemize}
\elemma

\prf (i) Existence: We define the $i$-th neighbor $S[i]$ as follows:
\begin{itemize}
\item[(a)] $S[0]=[z_1,\ldots,z_n,x_0]$, where $x_0=z_n+(z_1-z_0)$.
\item[(b)] $0<i<n$: Take $S[i]=[z_0,\ldots,z_{i-1},x_i,z_{i+1},\ldots,z_n]$, where $x_i=z_{i-1}+(z_{i+1}-z_i)$.
\item[(c)] $S[n]=[x_n,z_0,\ldots,z_{n-1}]$, where $x_n=z_0-(z_n-z_{n-1})$.
\end{itemize}
It is obvious that $\#(S\cap S[i])=n$ in all three cases. In the three cases, the distances between consecutive
points of $S[i]$ are given by (a) $e_{\alpha(2)},\ldots,e_{\alpha(n)},e_{\alpha(1)}$,
(c) $e_{\alpha(n)},e_{\alpha(1)},\ldots,e_{\alpha(n-1)}$, and (b) 
$e_{\alpha(1)},\ldots,e_{\alpha(i-1)},e_{\alpha(i+1)},e_{\alpha(i)},e_{\alpha(i+2)},\ldots,e_{\alpha(n)}$. (Figure
\ref{fig-kulpa} should make this quite clear.) Thus
$S[i]$ is a legal $n$-simplex for each $i\in\{0,\ldots,n\}$. Uniqueness: It remains to show that these are the
only ways of defining $S[i]$ consistently with $S\cap S[i]=\{z_0,\ldots,z_{i-1},z_{i+1},\ldots,z_n\}$. 
In the cases $i=0$ or $i=n$, the latter condition implies that $S[i]$ has a string of $n$ consecutive $z_i$'s
in common with $S$, and therefore also their differences given by $n-1$ mutually different vectors $e_j$. This
means that only one such vector is left, and the only way to use it so that $S[i]$ is an $n$-simplex different
from $S$ is to use 
it at the other end of the string of $z$'s. This shows the uniqueness of the above definitions in cases (a) and (c).
In the case $0<i<n$, $S$ and $S[i]$ have two corresponding substrings of $z$'s. A little thought shows that
the order of these two substrings must be the same in $S[i]$ as in $S$, so that all we can do is exchange two
adjacent difference vectors $e_{\alpha(i)}, e_{\alpha(i+1)}$, as done in the definition of $S[i]$ in case (b).

(ii) We must check whether $S[i]\subset C(k)$, which amounts to checking whether the new point $x_i$ is in
$I^n$. In case (a), we have $S=[z_1-e_{\alpha(1)},z_1,z_2,\ldots,z_n]$ and $S(0)=[z_1,\ldots,z_n,z_n+e_{\alpha(1)}]$.
If $F_0(s)=[z_1,\ldots,z_n]\subset C_j^\ve(k)$ then $z_1,\ldots,z_n$ all have the same $j$-coordinate $c$,
thus we must have $\alpha(1)=j$ and $c=1$ (since $S\subset I^n$). But then $S[0]\not\subset I^n$.
Conversely, if both $S$ and $S[0]$ are in $I^n$, then $z_1$ must have $\alpha(1)$-th coordinate $>0$ and $z_n$
must have $\alpha(1)$-th coordinate $<1$. All other coordinates of $z_1,\ldots,z_n$ are non-constant since the
vectors $e_{\alpha(2)},\ldots,e_{\alpha(n)}$ appear as differences. Thus $F_0(S)$ is not contained in any face
$I_i^\ve$. The cases (b) and (c) are checked similarly.
\qed

\noindent{\it Proof of Proposition \ref{prop-sperner}.}
For later use, we note the following fact (*): If $S\subset I^n$ satisfies 
$\varphi(S\cap I_i^\ve)=\{0,\ldots,n-1\}$ then $i=n$ and $\ve=-$. [The statement
$\varphi(S\cap I_i^\ve)=\{0,\ldots,n-1\}$ is contradicted by assumption (ii) if $\ve=+$  and 
  by (i) if $\ve=-$ and $i<n$.]

We call a subset $S\subset C(k)$ with $l+1$ elements \underline{full} if $\varphi(S)=\{0,\ldots,l\}$.
By (vi), a full $n$-simplex $S$ meets all $H_i^\pm$.  We will prove that the number $N_k$ of
full $n$-simplices in $C(k)$ is odd, thus non-zero, for all $k$.
The proof of $N_k\equiv 1\ (\mathrm{mod}\ 2)$ proceeds by induction over the dimension $n$ of $C(k)$ (for
fixed $k$). For $n=0$ we have $C(k)=\{0\}$, and there is exactly one full $n$-simplex, namely
$S=[z_0=0]$. Thus $N_0=1$. 

For an $n$-simplex $S\subset C(k)$, let $N(S)$ denote the number of full $(n-1)$-faces of $S$. If $S$ is
full then $N(S)=1$. [Since $\varphi(S)=\{0,\ldots,n\}$ and the only full $(n-1)$-face is obtained by
  omitting the unique $z_i$ for which $\varphi(z_i)=n$.] If $S$ is not full then $N(S)=0$ in the case
$\{0,\ldots,n-1\}\not\subset\varphi(S)$ [since omitting a $z_i$ cannot give a full (n-1)-face] or $N(S)=2$
in the case $\varphi(S)=\{0,\ldots,n-1\}$ [since there are $i\ne i'$ such that $z_i=z_{i'}$, so that $S$ becomes
  full upon omission of either $z_i$ or $z_{i'}$].
Thus 
\be N_k\equiv \sum_S N(S)\ \ (\mathrm{mod}\ 2), \label{eq1}\ee
where the summation extends over all $n$-simplices in $C(k)$.

Now by the Lemma, an $(n-1)$-face $F\subset C(k)$ belongs to one or two $n$-simplices in $C(k)$, depending on
whether $F\subset\del C(k)$ or not. Thus only the full faces $F\subset\del C(k)$ contribute to (\ref{eq1}): 
\[ N_k\equiv \# \{ F\subset\del C(k)\ \mathrm{full}\ (n-1)\mathrm{-face}\}\ \ (\mathrm{mod}\ 2). \]
If $F\subset\del C(k)$ is a full $(n-1)$-face then (*) implies $F\subset C_n^-(k)$. 
We can identify $C_n^-(k)=C(k)\cap I_n^-$ with $C_{n-1}(k)$, and under this identification $F$ is a
full $(n-1)$-simplex in $\7Z_k^{n-1}$. Thus $N_k\equiv N_{k-1}\ (\mathrm{mod}\ 2)$.
By the induction hypothesis, $N_{k-1}$ is odd, thus $N_k$ is odd.

Thus there is a full $n$-simplex $S=[z_0,\ldots,z_n]$, and if 
$C=\{ z_0+\sum_{i=1}^n a_i e_i\ | \ a\in\{0,1\}^n\}$ we have $S\subset C$ so that $\varphi(C)=\{0,\ldots,n\}$.
\qed

\brem Combinatorial proofs of Sperner's lemma, whether simplicial or cubical, may appear
mysterious. Homological proofs tend to be more transparent, cf.\ e.g.\ \cite{ivanov} in the simplicial case,
but here the point of course is to avoid the heavy machinery of homology.
\erem

\section{Higher connectedness of the cube}
Now we are in a position to prove, still following \cite{kulpa}, the following 
beautiful theorem, which for $n=1$ is just the connectedness of $[0,1]$:

\btheor \label{theor-kulpa} For $i=1,\ldots,n$, let $H_i^+, H_i^-\subset I^n$ be closed sets such that for
all $i$ one has $I_i^\pm\subset H_i^\pm$ and $H_i^-\cup  H_i^+=I^n$. Then $\bigcap_i(H_i^-\cap H_i^+)\ne\emptyset$.  
\etheor

\prf We define $F_0=I^n$ and $F_i=H_i^+\backslash I_i^-$ for all
$i\in\{1,\ldots,n\}$. Now define a map $\varphi:I^n\rarr\{0,\ldots,n\}$ by
\[ \varphi(x)=\max\left\{j:\ x\in\bigcap_{k=0}^j F_k\right\}. \]
Since $I_i^-\cap F_i=\emptyset$, we have $x\in I_i^-\ \impl\ \varphi(x)<i$. On the other hand,
if $x\in I_i^+$ then $\varphi(x)\ne i-1$. Namely, $\varphi(x)=i-1$ would mean that
  $x\in\bigcap_{k=0}^{i-1} F_k$ and $x\not\in F_i$. But this is impossible since $x\in I_i^+\subset H_i^+$ and 
  $I_i^+\cap I_i^-=\emptyset$, thus $x\in F_i$. 

By the above, the restriction of $\varphi$ to $C(k)\subset I^n$ satisfies the assumptions of Proposition
\ref{prop-sperner}. 
Thus for every $k\in\7N$ there is a subcube $C_k\subset C(k)$ such that $\varphi(C_k)=\{0,\ldots,n\}$.
Now, if  $\varphi(y)=i\in\{1,\ldots,n\}$ then $y\in F_i=H_i^+\backslash I_i^-\subset H_i^+$. 
On the other hand, if $\varphi(x)=i-1\in\{0,\ldots,n-1\}$ then $x\not\in F_i=H_i^+\backslash I_i^-$, which is
equivalent to $x\not\in H_i^+\vee x\in I_i^-$. The first alternative implies $x\in H_i^-$ 
  (since $H_i^+\cup H_i^-=I^n$), as does the second (since $I_i^-\subset H_i^-$). In either case, $x\in H_i^-$.
Combining these facts, we find that $\varphi(C_k)=\{0,\ldots,n\}$ implies
$C_k\cap H_i^+\ne\emptyset\ne C_k\cap H_i^-\ \forall i$, thus the subcube $C_k$ meets all the $H_i^\pm$. We
clearly have $\diam(C_k)=\sqrt{n}/k$, and since $k$ was arbitrary, the following lemma applied to 
$X=I^n,\ S_k=C_k$ and $\{ K_1,\ldots,K_{2n}\}=\{ H_i^\pm\}$ completes the proof.
\qed

\blemma \label{lem-inters} Let $(X,d)$ be a metric space and $\{K_i\subset X\}_{i\in I}$ compact subsets. Let
$\{ S_k\subset X\}_{k\in\7N}$ such that $\diam(S_k)\stackrel{k\rarr\infty}{\longrightarrow} 0$ and 
for all $k\in\7N,\ i\in I$ one has $S_k\cap K_i\ne\emptyset$. Then $\bigcap_i K_i\ne\emptyset$.
\elemma

\prf It is sufficient to prove the claim in the case where $I$ is finite. In the general case, this then
implies that the family $\{K_i\}_{i\in I}$ has the finite intersection property and another invocation of
compactness gives $\cap_i K_i\ne\emptyset$.
Thus let $\{K_1,\ldots,K_n\}$ be given and consider $K=\prod_i K_i$
equipped with the metric $d_K(x,y)=\sum_i d(x_i,y_i)$. For every $k\in\7N$ and
$i\in\{1,\ldots,n\}$, choose an $x_{k,i}\in S_k\cap K_i$ and define $x_k=(x_{k,1},\ldots,x_{k,n})\in K$. 
By compactness of $K$ there exists a point $z=(z_1,\ldots, z_n)\in K$ every neighborhood
of which contains $x_k$ for infinitely many $k$. Now, 
$d(z_i,z_j)\le d(z_i,x_{k,i})+d(x_{k,i},x_{k,j})+d(x_{k,j},z_j)\le 2d_K(z,x_k)+\mathrm{diam(S_k)}$. Since by
construction every neighborhood of $z$ contains 
points $x_k$ with arbitrarily large $k$, we can make both terms on the r.h.s.\ arbitrarily small and conclude
that $z=(x,\ldots,x)$ for some $x\in X$. Since $z_i\in K_i$ for all $i$, we have $x\in\cap_i K_i$, and are done. 
\qed

\section{The dimension of $I^n$} \label{sec-dim}
\bdefin \label{def-separate} If $A,B,C\subset X$ are closed sets such that $X\backslash C=U\cup V$, where
$U,V$ are disjoint open sets such that $A\subset U$ and $B\subset V$, we say that $C$ \underline{separates}
$A$ and $B$. 
\edefin

The following result plays an essential r\^ole in virtually all accounts of dimension theory. While it is
usually derived from Brouwer's fixed-point theorem, we obtain it more directly as an obvious corollary of Theorem 
\ref{theor-kulpa}.

\bcoro \label{coro-kulpa}
Whenever $C_1,\ldots,C_n\subset I^n$ are closed sets such  that $C_i$ separates $I_i^-$ and $I_i^+$ for each
$i$, then $\bigcap_i C_i\ne\emptyset$. 
\ecoro

\prf In view of Definition \ref{def-separate}, we have open sets $U_i^\pm$ such that $I_i^\pm\subset U_i^\pm$,  
$U_i^+\cap U_i^-=\emptyset$ and $U_i^+\cup U_i^-=X\backslash C_i$ for all $i$. Define 
$H_i^\pm=U_i^\pm\cup C_i$. Then $X\backslash H_i^\pm=U_i^\mp$, thus $H_i^\pm$ is closed. By construction,
$I_i^\pm\subset H_i^\pm$ and $H_i^+\cup H_i^-=I^n$, $H_i^+\cap H_i^-=C_i$, for all $i$. Now Theorem
\ref{theor-kulpa} gives $\cap_i C_i=\cap_i(H_i^-\cap H_i^+)\ne\emptyset$.
\qed

The preceding result will provide a lower bound on the dimension of $I^n$. The next result, taken from
\cite{mill}, will provide the upper bound:

\bprop \label{prop-mill} Let $A_1,B_1,\ldots, A_{n+1}, B_{n+1}\subset I^n$ be closed sets
such that $A_i\cap B_i=\emptyset$ for all $i$. Then for all $i$ there exist closed sets $C_i$ separating
$A_i$ and $B_i$  and satisfying $\bigcap_i C_i=\emptyset$.  
\eprop

\prf Pick real numbers $r_1, r_2,\ldots$ such that $r_i-r_j\not\in\7Q$ for $i\ne j$. 
(It suffices to take $r_k=k\sqrt{2}$.)
Then the sets $E_i=r_i+\7Q$ are mutually disjoint dense subsets of $\7R$.

Let $A,B\subset I^n$ be disjoint closed sets and $E\subset\7R$ dense. Then for every $x\in A$ we can find
an open neighborhood $U_x=I^n\cap\prod_{i=1}^n (a_i,b_i)$ with $a_i,b_i\in E$ such that $\ol{U_x}$ is disjoint
from $B$.  Since $A\subset I^n$ is closed, thus compact, there are $x_1,\ldots,x_k\in A$ such that 
$U=U_{x_1}\cup\cdots\cup U_{x_k}\supset A$. Now $C=\del U\subset I^n$ is closed and 
$X\backslash C=U\cup V$, where $V=I^n\backslash\ol{U}$. Now $U,V$ are open and disjoint such that $A\subset U,
B\subset V$, thus $C$ separates $A$ and $B$. If $x\in\del(I^n\cap\prod_i(a_i,b_i))$ then at least one of the
coordinates $x_i$ of $x$ equals $a_i$ or $b_i$, and thus is in $E$. Now,  
$C=\del U\subset\del U_{x_1}\cup\cdots\cup \del U_{x_k}\subset\{ x\in I^n\ | \ \exists j\in\{1,\ldots,n\}: x_j\in E\}$. 

We can thus find, for each pair $(A_i,B_i)$ a closed set 
$C_i\subset\{ x\in I^n\ | \ \exists j: x_j\in E_i\}$ that separates $A_i$ and $B_i$. Let now $x\in\cap_i C_i$.  
Then for every $i\in\{1,\ldots,n+1\}$ there is a $j_i\in\{1,\ldots,n\}$ such that $x_{j_i}\in E_i$. By the
pigeonhole principle there are $i,i'\in\{1,\ldots,n+1\}$ such that $i\ne i'$ and $j_i=j_{i'}=j$. But this
means that  $x_j\in E_i\cap E_{i'}\in\emptyset$, which is absurd. Thus $\cap_i C_i=\emptyset$.
\qed

Proposition \ref{prop-mill} should be compared with Corollary \ref{coro-kulpa}. In order to do this
systematically, the following is convenient:

\bdefin \label{sdim} Let $X$ be a topological space. We define the \underline{separation-dimension} 
$s$-$\dim(X)\in\{-1,0,1,\ldots,\infty\}$ as follows:
\begin{itemize}
\item We put $s$-$\dim(X)=-1$ if and only if $X=\emptyset$.
\item If $X\ne\emptyset$ and $n\in\7N_0$, we say that $s$-$\dim(X)\le n$ if, given closed sets
  $A_1,B_1,\ldots, A_{n+1}, B_{n+1}$ such that $A_i\cap B_i=\emptyset\ \forall i$, there exist closed
  $C_i$ separating $A_i$ and $B_i$ and satisfying $\cap_i C_i=\emptyset$. 
(This is consistent: If $s$-$\dim(X)\le n$ and $n<m$ then $s$-$\dim(X)\le m$.)
\item If $s$-$\dim(X)\le n$ holds, but $s$-$\dim(X)\le n-1$ does not, we say $s$-$\dim(X)=n$.
\item If there is no $n\in\7N$ such that $s$-$\dim(X)\le n$ then $s$-$\dim(X)=\infty$.
\end{itemize}
\edefin

\blemma If $X$ and $Y$ are homeomorphic then $s$-$\dim(X)=s$-$\dim(Y)$. \elemma

\prf Obvious since $s$-$\dim(X)$ is defined in terms of the topology of $X$. \qed

\btheor \label{theor-dimI}$s$-$\dim(I^n)=n$.
\etheor

\prf Proposition \ref{prop-mill} implies $s$-$\dim(I^n)\le n$. On the other hand, it is clear that
$s$-$\dim(I^n)\ge n$ holds if and only if there are closed sets $A_1,B_1,\ldots,A_n,B_n\subset X$ satisfying 
$A_i\cap B_i=\emptyset\ \forall i$ such that any closed sets $C_i$ separating $A_i$ and $B_i$ satisfy
$\cap_iC_i\ne\emptyset$. This is exactly what is asserted by Corollary \ref{coro-kulpa}.
\qed

\bcoro We have $I^n\cong I^m$ if and only if $n=m$.
\ecoro

\brem \label{rem-hist} 1. A family $\{ (A_1,B_1),\ldots,(A_{n+1},B_{n+1})\}$ as above is called
`essential'. Therefore one could also speak of the `essential family dimension', but this does not seem to be
in widespread use.

2. In 1938, Eilenberg and Otto \cite{EO} proved that the separation dimension coincides with the covering
dimension in the case of separable metrizable spaces. This was generalized to normal spaces by Hemmingsen
(1946). For a modern proof and more history see \cite{engelking3}, in particular Theorem 3.2.6.

3. It is a triviality to prove $s$-$\dim(Y)\le s$-$\dim(X)$ for closed $Y\subset X$, implying
$s$-$\dim(\7R^n)\ge n$ and $s$-$\dim(S^n)\ge n$. In order to prove the converse inequalities, and thereby the
invariance of dimension for spheres and Euclidean spaces, one needs a `sum theorem' for the separation
dimension. This is the statement that if $X=\cup_{i\in \7N}Y_i$ with $Y_i\subset X$ closed and
$s$-$\dim(Y_i)\le d$ then $s$-$\dim(X)\le d$. Such a result follows from the combination of the 
first half of this remark and the known sum theorem for the covering dimension, cf.\ 
\cite[Theorem 3.1.8]{engelking3}. (This is ``just'' point set topology, with no simplicial or combinatorial 
methods involved.) However, from an aesthetic perspective it would be desirable to give a direct proof of the
sum theorem for $s$-$\dim$. 
\erem

\appendix
\section{The theorems of Poincar\'e-Miranda and Brouwer} \label{ss-brouwer}
In order to further illustrate the power of Theorem \ref{theor-kulpa}, and for the benefit of the reader, we 
include deductions of many important classical results about continuous functions on $I^n$, following
\cite{kulpa}. But we emphasize that none of this is needed for the treatment of the invariance of dimension in 
Section \ref{sec-dim}.

\bcoro [Poincar\'e-Miranda theorem] \label{coro-miranda}
Let $f=(f_1,\ldots,f_n)\in C(I^n,\7R^n)$. If
$f_i(I_i^-)\subset(-\infty,0]$, $f_i(I_i^+)\subset[0,\infty)$ for all $i$, then there is $x\in I^n$ such that
    $f(x)=0$. 
\ecoro

\prf Put $H_i^-=f_i^{-1}((-\infty,0]), H_i^+=f_i^{-1}([0,\infty))$. Then clearly 
$I_i^\pm\subset H_i^\pm$ and $H_i^-\cup  H_i^+=I^n$, for all $i$. By Theorem \ref{theor-kulpa}, there exists
$x\in\cap_i(H_i^-\cap H_i^+)$, and it is clear that $f(x)=0$.
\qed

\bcoro [Brouwer's fixed-point theorem] \label{coro-brouwer} Let $g\in C(I^n,I^n)$. Then there exists $x\in I^n$
such that $g(x)=x$. (I.e., $I^n$ has the fixed-point property.)
\ecoro

\prf Put $f(x)=x-g(x)$. Then the assumptions of Corollary \ref{coro-miranda} are satisfied, so that there is
$x\in I^n$ for which $f(x)=0$. Thus $g(x)=x$.
\qed

\bcoro \label{coro-surj} Let $g\in C(I^n,I^n)$ satisfy $g(I_i^\pm)\subset I_i^\pm\ \forall
i$. (E.g.\ $g\restr\del I^n=\id$.) Then $g(I^n)=I^n$.
\ecoro

\prf Let $p\in I^n$, and put $f(x)=g(x)-p$.  Then $f$ satisfies the assumptions of Corollary
\ref{coro-miranda}, thus there is $x\in I^n$ with $f(x)=0$. This means $g(x)=p$, so that $g$ is surjective.
\qed

\brem 1. The history of the above results is quite convoluted and interesting. See the introduction of
\cite{kulpa} for a glimpse.

2. Since compact convex subsets of $\7R^n$ are homeomorphic to $D^m$ for some $m\le n$, they also have the
fixpoint property. The convexity assumption cannot be omitted as is shown by a non-trivial rotation of
$S^1\subset\7R^2$. On the other hand, Corollary \ref{coro-surj} and the resulting non-existence of retractions
to the boundary extend to arbitrary compact subsets of $\7R^n$, cf.\ \cite{kulpa}.

3. The Poincar\'e-Miranda theorem seems to be much more popular with analysts than with topologists. One may
indeed argue that Brouwer's theorem is more fundamental, asserting the fixed-point property of any $n$-cell 
irrespective of its shape (e.g.\ disk, cube or simplex). But apart from being a particularly natural higher
dimensional generalization of the intermediate value theorem, the Poincar\'e-Miranda theorem often is the more 
convenient point of departure for other proofs. (The Poincar\'e-Miranda theorem can be deduced from Brouwer's
theorem, cf.\ \cite{miranda}, but the argument is more involved. Cf.\ also \cite[p.\ 118]{ist}.) 

4. It is also true that the Poincar\'e-Miranda theorem implies Theorem \ref{theor-kulpa}: Given $H_i^\pm$ as
in the latter, the functions $f_i(x)=\dist(x,H_i^-)-\dist(x,H_i^+)$ are continuous and satisfy
$f_i(I_i^-)\subset[-1,0]$, $f_i(I_i^+)\subset[0,1]$. Now the Poincar\'e-Miranda theorem gives an $x\in I^n$
such that $f(x)=0$. The assumption $H_i^-\cup H_i^+=I^n$ implies that $\dist(x,H_i^-), \dist(x,H_i^+)$ cannot
both be non-zero. Thus $f_i(x)=0$ is equivalent to $x\in H_i^-\cap H_i^+$. Thus
$x\in\cap_i(H_i^-\cap H_i^+)$. 

5. Combining the above facts and some other well-known implications, we see that the following statements are
`equivalent' (in the sense of being easily deducible from each other): 
\begin{itemize}
\item[(i)] the non-existence of a retraction $r: D^n\rarr\del D^n$, 
\item[(ii)] the non-contractibility of $\del D^n =S^{n-1}$, 
\item[(iii)] $\pi_{n-1}(S^{n-1})\ne 0$,
\item[(iv)] the fixed point property of $D^n$, 
\item[(v)] the statement $f:D^n\rarr D^n,\ f\restr\del D^n=\id\ \impl\ f(D^n)=D^n$, 
\item[(vi)] the Poincar\'e-Miranda theorem, 
\item[(vii)] Theorem \ref{theor-kulpa}. 
\end{itemize}
(A similar statement appears in \cite[Theorem 6.6.1]{TD}, where however Corollary \ref{coro-kulpa} is listed 
instead of the more convenient Theorem \ref{theor-kulpa} and correspondingly, the Poincar\'e-Miranda type
theorem given there makes the stronger assumptions $f_i(I_i^-)\subset(-\infty,0)$,
$f_i(I_i^+)\subset(0,\infty)\ \forall i$.)

We observe that (i)-(vi) all involve continuous maps, whereas (vii) only involves the faces and the topology
of the cube $I^n$. The latter therefore seems closest in spirit to point set topology, as is also supported by
its interpretation as higher-dimensional connectedness and its r\^ole in the above proof of dimension invariance. 
\erem

\end{document}